\magnification=\magstep1 \baselineskip=1.3\baselineskip

\def\vphi{\varphi}
\def\eps{\varepsilon}
\def\R{{\bf R}}
\def\prt{{\partial}}
\def\df{{\mathop {\ =\ }\limits^{\rm{df}}}}
\def\ol{\overline}
\def\dist{{\mathop {\rm dist}}}
\def\diam{{\mathop {\rm diam}}}
\def\wt{\widetilde}
\def\wh{\widehat}
\def\bone{\hbox{\bf 1}}
\def\n{{\bf n}}

\def\qed{{\hfill $\square$ \bigskip}}
\def\sqr#1#2{{\vcenter{\vbox{\hrule height.#2pt
        \hbox{\vrule width.#2pt height#1pt \kern#1pt
           \vrule width.#2pt}
        \hrule height.#2pt}}}}
\def\square{\mathchoice\sqr56\sqr56\sqr{2.1}3\sqr{1.5}3}

\input epsf.tex
\newdimen\epsfxsize
\newdimen\epsfysize

\centerline{\bf THE ``HOT SPOTS'' PROBLEM }
 \centerline{\bf IN PLANAR DOMAINS WITH ONE HOLE}

 \footnote{$\empty$}{\rm Research partially supported by
NSF grant DMS-0303310.}
\bigskip

\vskip0.4truein \centerline{\bf Krzysztof Burdzy}

\vskip1truein

{\narrower \noindent
 {\bf Abstract}. There exists a planar domain with piecewise
 smooth boundary and one hole such that the second
 eigenfunction for the Laplacian with Neumann boundary
 conditions attains its maximum and minimum inside the domain.

 }

\vskip1truein

\noindent {\bf 1. Introduction}.

We will be concerned with bounded planar domains with piecewise
smooth boundaries. The Laplacian with Neumann boundary conditions
in such a domain has a discrete spectrum (see, e.g., [BB1]).
Recall that the first eigenvalue is equal to 0 and let $\lambda$
denote the second eigenvalue. Our main result is the following.

\bigskip
\noindent{\bf Theorem 1.1}. {\sl There exists a planar domain $D$
with one hole, such that the second Neumann eigenvalue is simple
(i.e.,  the subspace of $L^2$ corresponding to $\lambda$ is
one-dimensional) and the corresponding eigenfunction $\vphi$
satisfies
 $$\inf_{x\in D} \vphi(x) < \inf_{x\in \prt D} \vphi(x)
 \leq \sup_{x\in \prt D} \vphi(x) < \sup_{x\in D} \vphi(x).
 \eqno(1.1)
 $$
}

\bigskip

The ``hot spots'' conjecture of J.~Rauch, proposed in 1974,
states, roughly speaking, that the second Neumann eigenfunction
attains its maximum on the boundary of a Euclidean domain. The
conjecture is false at this level of generality (see [BB2] and
[BW]) but it is true for some classes of domains (see [A], [AB],
[BB1], [JN], [K], [P]). The counterexample given in [BW] is a
planar domain with two holes and it suggests, in the intuitive
sense, that any planar domain where the hot spots conjecture fails
must have at least two holes. Theorem 1.1 shows that this is not
true. Theorem 1.1 is also a small step towards understanding of
the ``hot spots'' problem for domains with no holes. The first
part of the following version of the ``hot spots'' conjecture was
stated by Kawohl [K] while the second part is our own.

\bigskip
\noindent{\bf Conjecture 1.2}. {\sl (i) The second Neumann
eigenfunction attains its maximum on the boundary if $D$ is any
convex domain in $\R^n$, for any $n\geq 1$.

(ii) The second Neumann eigenfunction attains its maximum on the
boundary if $D$ is a simply connected planar domain. }

\bigskip

Counterexamples to the ``hot spots'' conjecture presented in [BW]
and [BB2] involved domains with bizarre shapes (the shapes were
unusual for technical reasons). The counterexample given in this
paper is rather simple (see Fig.~1 in the next section) so it
shows that the hot spots conjecture fails in some ``ordinary''
domains.

It was pointed out in [BB2] that it would be rather easy to
construct a two-dimensional manifold with a boundary (see
Figs.~2.1 and 2.2 in [BB2]) based on the same idea as that in
[BW], with the property that both maximum and minimum of the
second Neumann eigenfunction lie inside the manifold. The example
given in this article is much harder, from the intuitive point of
view, because a similar distortion of the example (i.e., a
two-dimensional manifold of a similar shape) would not be any
easier to deal with than the planar domain itself.

One of the goals of this paper is to develop new techniques for
studying the ``hot spots'' problem. Many of the articles cited
above converted the ``hot spots'' problem for eigenfunctions with
Neumann boundary conditions to a mixed boundary problem, by
cutting the domain into two subdomains along the nodal line (i.e.,
zero line) for the second Neumann eigenfunction (the nodal line
becomes a part of the boundary with the Dirichlet boundary
conditions). This technique proved to be very fruitful and we will
apply it in this paper. However, when the geometry of the domain
is not very simple, it is either hard to find the location of the
nodal line or to incorporate the nodal line into the argument. The
main part of the proof of Theorem 1.1 will be based on cutting the
domain along a level line of the second eigenfunction. This
modification makes it necessary to develop arguments that are more
quantitative than qualitative in nature, as compared to the
existing proofs. Of course, our proofs will include many ideas
from the existing literature, for example, [BB1] and [BW].

We will now briefly describe the idea of the proof of Theorem 1.1.
The domain $D$ depicted in Fig.~1 (see the next section) has two
axes of symmetry. First we will show that $\vphi$ is symmetric
with respect to one of them and antisymmetric with respect to the
other one. Hence, it is enough to analyze the upper right quarter
of the domain; let us call this subdomain $D_1$. The set $D_1$ is
a very thin ``tube'' of slightly variable width. The point $(0,0)$
lies on the boundary of $D_1$ and it is enough to show that
$\vphi$ is strictly larger at $(0,0)$ than at any point in $\prt D
\cap \prt D_1$. The point $(0,0)$ is the most distant point from
the other end of the tube $D_1$, in the sense that a reflected
Brownian motion in $D_1$ starting from $(0,0)$ will take longer
(on average) to reach the other end of $D_1$ than a reflected
Brownian motion starting from any other point of $\prt D \cap \prt
D_1$. This probabilistic statement can be translated into an
estimate needed for the proof of (1.1).

\bigskip
\noindent{\bf 2. Proofs}.

Our proofs will rely to large extent on techniques developed in
[BB1] and other papers. We will be brief at many places to keep
this article short. We ask the reader to consult [BB1] and other
articles cited below for more details.

An open disc with center $x$ and radius $r$ will be denoted
$B(x,r)$. We will identify points $x\in\R^2$ with vectors
$\overrightarrow{(0,0),x}$ and complex numbers $x=r e^{i\theta}$.
The angle between $x=r_x e^{i\theta_x}$ and $y=r_y e^{i\theta_y}$,
i.e., $\theta_x-\theta_y$, will be denoted $\angle (x,y)$. We will
write $\angle(x)$ instead of $\angle(x, (1,0))$, i.e., $\angle(x)$
will denote the angle formed by the vector $x$ with the positive
horizontal semi-axis. We will use the convention that
$\angle(x,y)\in(-\pi, \pi]$. For any process $Z_t$ we will denote
the hitting time of a set $A$ by $T^Z_A$, i.e., $T^Z_A =
\inf\{t\geq 0: Z_t \in A\}$. The superscript will be dropped if no
confusion may arise.

Our definition of a domain $D\subset \R^2$ satisfying (1.1) will
involve a parameter $\eps \in(0,1/4)$. The value of $\eps$ will be
chosen later and should be thought of as a very small number; it
will be suppressed in the notation. Let $A_1 $ be a convex
polygonal domain with the consecutive vertices $(0,-\eps),
(0,\eps), (1, 2\eps),(2, \eps_0),(2,-\eps_0)$ and $(1,-2\eps)$,
where $\eps_0 \in (0,\eps)$. The value of the parameter $\eps_0$
will be specified later. Let $C_1$ be a polygonal Jordan arc
inside
 $$(B((2,-1), 1+2\eps_0) \setminus B((2,-1), 1+\eps_0/2))
 \cap\{(x_1,x_2): x_1\geq 2, x_2 \geq -1\},$$
with endpoints $(2,\eps_0)$ and $(3+\eps_0, -1)$, and such that
for any line segments $\ol{x,y}, \ol{y,z}\subset C_1$ we have
$|\angle (y-x,z-y )| \leq \eps_0$. Similarly, let $C_2$ be a
polygonal Jordan arc inside
 $$(B((2,-1), 1-\eps_0/2) \setminus B((2,-1), 1-2\eps_0))
 \cap\{(x_1,x_2): x_1\geq 2, x_2 \geq -1\},$$
with endpoints $(2,-\eps_0)$ and $(3-\eps_0, -1)$, and such that
for any line segments $\ol{x,y}, \ol{y,z}\subset C_2$ we have
$|\angle (y-x,z-y )| \leq \eps_0$. Let $A_2$ be an open domain
whose boundary consists of $C_1$, $C_2$, and line segments
$\ol{(2,\eps_0), (2,-\eps_0)}$ and $\ol{(3+\eps_0, -1),(3-\eps_0,
-1)}$. Let $A_3 = A_1 \cup A_2$, let $A_4$ be the symmetric image
of $A_3$ with respect to the line $\{(x_1,x_2): x_2=-1\}$, and let
$A_5$ and $A_6$ be the symmetric images of $A_3$ and $A_4$ with
respect to $\{(x_1,x_2): x_1=0\}$. Finally we let $D$ be the
interior of the closure of $A_3 \cup A_4 \cup A_5 \cup A_6$. A
schematic drawing of $D$ is presented in Fig.~1. The polygonal
lines $C_1$ and $C_2$ are very close to circular arcs so they are
represented graphically as such. A substantial part of the
argument will be focused on a subdomain $D_1$ of $D$ depicted in
Fig.~1.

\bigskip
\vbox{ \epsfxsize=3.5in
  \centerline{\epsffile{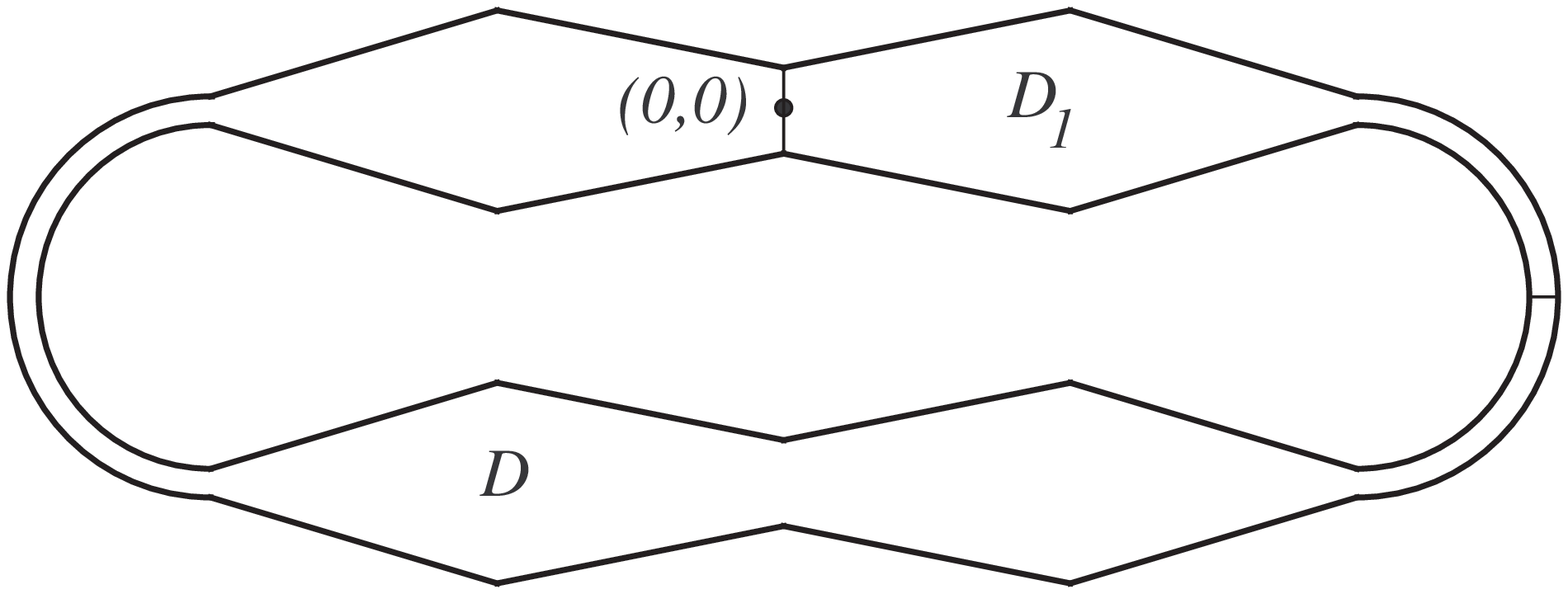}}

\centerline{Figure 1. (Drawing not to scale.)} }
\bigskip

We will now review a few basic facts about reflected Brownian
motion and ``synchronous'' couplings. Let $\n(x)$ denote the unit
inward normal vector at $x\in\prt D$. Let $W$ be standard planar
Brownian motion, $x,y\in \ol D$, and consider the following
Skorohod equations,
 $$\eqalignno{
 X_t &= x + W_t + \int_0^t  \n(X_s) dL^X_s, &(2.1)\cr
 Y_t &= y + W_t + \int_0^t  \n(Y_s) dL^Y_s. &(2.2) }$$
Here $L^X$ is the local time of $X_t$ on $\prt D$, i.e., a
non-decreasing continuous process which does not increase when $X$
is in $D$: $\int_0^\infty \bone_{D}(X_t) dL^X_t = 0$, a.s.
Equation (2.1) has a unique pathwise solution $(X_t,L^X_t)$ such
that $X_t \in \ol D$ for all $t\geq 0$ (see [LS]). The ``reflected
Brownian motion'' $X$ is a strong Markov process. The same remarks
apply to (2.2), so $(X, Y)$ is also strong Markov. We will call
$(X, Y)$ a ``synchronous coupling.'' Note that on any interval
$(s,t)$ such that $X_u \in D$ and $Y_u \in D$ for all $u \in
(s,t)$, we have $X_u - Y_u = X_s - Y_s$ for all $u \in (s,t)$.

\bigskip

Recall that $\lambda$ denotes the second Neumann eigenfunction in
$D$.

\bigskip

\noindent{\bf Lemma 2.1}. {\sl For any $c_1>0$ there exists
$c_2>0$ such that if $\eps_0\leq c_2\eps$ then $\lambda\leq c_1$
and $\lambda$ is simple.

}

\bigskip
\noindent{\bf Proof}. Let $r = \eps_0/(2 \eps - \eps_0)$ and note
that the point $y \df (2+r, 0)$ lies at the intersection of
straight lines passing through the line segments $\ol{(1, 2\eps),
(2, \eps_0)}$ and $\ol{(1, -2\eps), (2, -\eps_0)}$. Let
$K=\{x=(x_1,x_2)\in D: |x_1| \leq 2, x_2>-1\}$,  $K_1 =
B(y,2r)\cap K$, $K_2= \prt B(y,1/2)\cap K$ and $K_3= \prt
B(y,1)\cap K$. Let $X$ be a reflected Brownian motion in $D$ with
$X_0 \in K_2$. Let $T_0=0$, and for $k\geq 1$ let
 $$\eqalign{
 S_k &= \inf\{t\geq T_{k-1}: X_t \in K_1 \cup K_3\},\cr
 T_k & = \inf\{t \geq S_k: X_t \in K_2\}.
 }$$
Let $R_t= \dist(X_t, y)$ and note that if $X$ is between $K_1$ and
$K_3$, the process $R$ is a 2-dimensional Bessel process because
the normal reflection of $X$ on $\prt D$ has no effect on $R$. It
follows that for any $p_1<1$, there exists $r_0>0$ so small that
if $r\leq r_0$ then
 $$P(X_{S_k} \in K_3 \mid {\cal F}_{T_{k-1}})
 = {\log (1/2) - \log (2r) \over
 \log 1 - \log (2r)}\geq p_1.$$
Moreover, for some $t_0>0$ not depending on $r$, $P(X_{S_k} \in
K_3, S_k-T_{k-1}>t_0 \mid {\cal F}_{T_{k-1}}) \geq p_1$.

Let $z = (-2-r, 0)$, $K'_1 = (B(y,2r)\cup B(z,2r))\cap K$, $K'_2=
(\prt B(y,1/2)\cup \prt B(z,1/2))\cap K$, $K'_3= (\prt B(y,1)\cup
\prt B(z,1))\cap K$, $T'_0=0$, and for $k\geq 1$ let
 $$\eqalign{
 S'_k &= \inf\{t\geq T'_{k-1}: X_t \in K'_1 \cup K'_3\},\cr
 T'_k & = \inf\{t \geq S'_k: X_t \in K'_2\}.
 }$$
By symmetry, $P(X_{S'_k} \in K'_3, S'_k-T'_{k-1}>t_0 \mid {\cal
F}_{T'_{k-1}}) \geq p_1$. By the repeated use of the strong Markov
property,
 $$P(T^X_{K'_1} \geq kt_0 \mid X_0 \in K'_2)
 \geq P(\bigcap_{1\leq j\leq k}\{X_{S'_j} \in K'_3,
 S'_j-T'_{j-1}>t_0\}\mid X_0 \in K'_2)\geq p_1^k.$$

Let $D_- =\{x=(x_1,x_2)\in D: x_2 < -1\}$ and $A = \prt D_- \cap
D$. Let $u(t,x)$ be the heat equation solution in $D$ with the
Neumann boundary conditions and the initial condition $u(0,x) = 1$
for $x\in D_-$ and $u(0,x) = 0$ otherwise. Note that $u$ can be
represented probabilistically as $u(t,x) = P(X_t \in D_- \mid X_0
= x)$. By the strong Markov property applied at $T^X_{A}$ and
symmetry, $P(X_t \in D_- \mid T^X_{A} < t) =1/2$, so for $x\in
K'_2$ and large $t$,
 $$\eqalign{
 u(t,x) &= P(X_t \in D_- \mid X_0 = x)\cr
 & = P(X_t \in D_- , T^X_{A} < t \mid X_0 = x)\cr
 & = (1/2) P( T^X_{A} < t \mid X_0 = x)\cr
 & = 1/2 - (1/2)P(T^X_{A} \geq t \mid X_0 =x)\cr
 & \leq 1/2 - (1/2)P(T^X_{K'_1} \geq t \mid X_0 =x )\cr
 & \leq 1/2 - (1/2)p_1^{t/(2t_0)} = 1/2 - (1/2)e^{(\log p_1/(2t_0)) t}.
 }$$
Since $p_1$ can be made arbitrarily close to $1$ by making $r$
small, $-\log p_1/(2t_0)>0$ can be arbitrarily close to $0$. By
symmetry, $u(t,x)$ converges to $1/2$ as $t\to\infty$. By
Proposition 2.1 of [BB1], $\sup_{x\in D}|u(t,x)-1/2|\leq c_3
e^{-\lambda t}$ for large $t$. Hence, $\lambda \leq -\log
p_1/(2t_0)$ and we see that for any $c_1>0$ we have $\lambda \leq
c_1$, provided $r$ is sufficiently small. If $c_2 < 1$ and
$\eps_0\leq c_2\eps$ then $r = \eps_0/(2 \eps - \eps_0)\leq c_2 $,
so $\lambda \leq c_1$ if we assume that $c_2$ is small. This
proves the first claim of the lemma.

The assertion that $\lambda $ is simple is totally analogous to
the claims proved in Sections 4 and 5 of [BW]. The proofs in [BW]
are based on the fact that the domain has a bottleneck and they
extend easily to our domain $D$. We leave the details to the
reader.
 \qed
\bigskip

We will assume from now on that $\eps_0$ and $ \eps $ are such
that $\lambda < 1$ and $\lambda $ is simple. Recall that the
``nodal line'' is the set of points $x$ such that $\vphi(x)=0$. We
will use the phrase ``nodal line'' even if the set of zeros of
$\vphi$ is not connected.

\bigskip

\noindent{\bf Lemma 2.2}. {\sl For any $\eps>0$ there is
$\eps_1\in(0,\eps)$ such that if $\eps_0\in (0,\eps_1)$ then the
following is true. The eigenfunction $\vphi$ is symmetric with
respect to the vertical axis and antisymmetric with respect to the
line $\{(x_1,x_2): x_2=-1\}$, i.e., for any $(x_1,-1+x_2)\in D$,
we have $\vphi(x_1,-1+x_2) = \vphi(-x_1, -1+x_2) = - \vphi(x_1,
-1-x_2) = - \vphi(-x_1, -1-x_2)$. It follows that the nodal line
is $\{(x_1,x_2)\in D: x_2 = -1\}$.

}

\bigskip
\noindent{\bf Proof}. The function $\vphi_1(x_1,x_2) \df
\vphi(x_1,x_2) + \vphi(-x_1,x_2)$ is an eigenfunction
corresponding to $\lambda$. If $\vphi_1$ is identically equal to
zero then $\vphi$ is antisymmetric with respect to the vertical
axis. If $\vphi_1$ is not identically equal to zero then it is a
constant multiple of $\vphi$ (because $\lambda$ is simple) and it
follows that $\vphi$ is symmetric with respect to the vertical
axis. A similar argument shows that either $\vphi$ is
antisymmetric with respect to $\{(x_1,x_2): x_2=-1\}$ or it is
symmetric with respect to this line.

An argument similar to that in Section 4 of [BW] shows that for
any fixed $\eps$, the nodal line cannot intersect the set
$\{(x_1,x_2)\in D: |x_1| \leq 1\}$ if  $\eps_0$ is sufficiently
small. By the Courant Nodal Line Theorem ([CH]), the nodal line
divides $D$ into two connected components. These facts taken
together with the symmetries described in the first paragraph of
the proof imply that the nodal line must be equal to
$\{(x_1,x_2)\in D: x_2 = -1\}$. It follows that $\vphi$ is
antisymmetric with respect to $\{(x_1,x_2): x_2=-1\}$ and it is
symmetric with respect to the vertical axis. \qed

\bigskip

The next lemma is a prelude to a theorem on geometric properties
of ``mirror'' couplings, to be defined later. The lemma is
concerned with convergence of a sequence of processes to the
reflected Brownian motion---we start with the construction of this
sequence. Suppose that $W$ is a planar Brownian motion, $x$ is a
point in the upper half-plane $D_* \df \{(y_1,y_2)\in\R^2:
y_2>0\}$, and $c_1<\infty$ is a constant. For every fixed
$\delta>0$, we will construct a process $X^\delta$ inductively.
Let $X^{\delta,1}$ be the reflected Brownian motion in $D_*$,
starting from $x$ and driven by $W$, in the sense of (2.1). Let
$T_0=0$, and $T_1\geq T_0$ be a stopping time such that
$X^{\delta,1}_{T_1}\in \prt D_*$ a.s. Let $V_1$ be a random
variable satisfying $|V_1| \leq c_1 \delta^2$, a.s. For the
induction step, suppose that the process $X^{\delta,j}$ is
defined, $T_j$ is a stopping time for $X^{\delta,j}$ such that
$T_j\geq T_{j-1}$, $X^{\delta,j}_{T_j}\in \prt D_*$ a.s., and
$V_j$ is a random variable satisfying $|V_j| \leq c_1 \delta^2$,
a.s. We define $\{X^{\delta, j+1}_t, t\geq T_j\}$ as the reflected
Brownian motion driven by $\{W_t, t\geq T_j\}$, starting at
$X^{\delta,j}_{T_j} + (V_j, \delta)$. Then we choose any
$X^{\delta, j+1}$-stopping time $T_{j+1}\geq T_j$ such that
$X^{\delta,j+1}_{T_{j+1}}\in \prt D_*$ a.s., and a random variable
$V_{j+1}$ satisfying $|V_{j+1}| \leq c_1 \delta^2$, a.s. The
process $X^\delta$ is defined by $X^\delta_t = X^{\delta, j}_t$
for $t\in[T_{j-1}, T_j)$. It is elementary to see that $T_j \to
\infty$ a.s., so $X^\delta_t$ is well defined for all $t\geq 0$
a.s.

\bigskip

\noindent{\bf Lemma 2.3}. {\sl The processes $X^\delta$ converge
in distribution to the reflected Brownian motion in $D_*$ as
$\delta \to 0$.

}

\bigskip
\noindent{\bf Proof}. Let us denote coordinates of processes as
follows, $W=(\wt W, \wh W)$ and $X^\delta = (\wt X^\delta, \wh
X^\delta)$. Note that $\wh X^\delta_t = \wh W _t + L^\delta_t$,
where $L^\delta$ is a non-decreasing process. It is elementary to
prove that $L^\delta$ converge to a process $L$ as $\delta\to0$,
on every time interval $[0,t_0]$, and the process $L$ is
non-decreasing, continuous and does not increase when $\wh W _t +
L_t>0$. By the uniqueness of the Skorohod decomposition, $\wh W _t
+ L_t$ is a one-dimensional reflected Brownian motion. Hence, $\wh
X^\delta$ converge to the Shorohod transform (in the sense of
(2.1)) of $\wh W$.

Fix a time interval $[0,t_0]$ and let $N=N(\delta)$ be the number
of $j$ with $T_j\leq t_0$. Note that $L^\delta_{t_0} \geq N
\delta$ and the jumps of $L^\delta$ occur only when $X^\delta$
approaches $\prt D_*$. Since $\{W_t, t\in[0,t_0]\}$ has a bounded
diameter a.s., it follows that there exists a random variable
$N_0<\infty$ such that $N\leq N_0/\delta$ for every
$\delta\in(0,1)$, a.s. This implies that that $\sum_{j\leq N}
|V_j| \leq c_1 N_0 \delta$. Since this random quantity converges
to $0$ in distribution, as $\delta\to0$, $\wt X^\delta$ converge
to $\wt W$.
 \qed

\bigskip

We will now review some properties of ``mirror couplings'' for
reflected Brownian motions which are relevant to our arguments.
These aspects of mirror couplings were originally developed in
[BK] and later applied in [BB1] and [BB2]. Our review is borrowed
from [BB2].

We start with the mirror coupling of two Brownian motions in
$\R^2$. Suppose that $x,y\in \R^2$ are symmetric with respect to a
line $M$. Let $X$ be a Brownian motion starting from $x$ and let
$Y_t$ be the mirror image of $X_t$ with respect to $M$ for $t \leq
T^X_M$. We let $Y_t = X_t$ for $t>T^X_M$. The process $Y$ is a
Brownian motion starting from $y$. The pair $(X,Y)$ is a ``mirror
coupling'' of Brownian motions.

Next we turn to the mirror coupling of reflected Brownian motions
in a half-plane $D_*$, starting from $x, y \in D_*$. Let $M$ be
the line of symmetry for $x$ and $y$. The case when $M$ is
parallel to $\prt D_*$ is essentially a one-dimensional problem,
so we focus on the case when $M$ intersects $\prt D_*$. By
performing rotation and translation, if necessary, we may suppose
that $D_*$ is the upper half-plane and $M$ passes through the
origin. We will write $x = (r^x, \theta^x)$ and $y = (r^y,
\theta^y)$ in polar coordinates. The points $x$ and $y$ are at the
same distance from the origin so $r^x = r^y$. Suppose without loss
of generality that $\theta^x < \theta^y$. We first generate a
2-dimensional Bessel process $R_t$ starting from $r^x$. Then we
generate two coupled one-dimensional processes on the
``half-circle'' as follows. Let $\wt \Theta^x_t$ be a
1-dimensional Brownian motion starting from $\theta^x$. Let $\wt
\Theta^y_t = - \wt \Theta^x_t +\theta^x +\theta^y$. Let
$\Theta^x_t$ be reflected Brownian motion on $[0,\pi]$,
constructed from $\wt \Theta^x_t$ by the means of the Skorokhod
equation.  Thus $\Theta_t^x$ solves the stochastic differential
equation $d\Theta_t^x=d\wt \Theta_t^x+dL_t$, where $L_t$ is a
continuous process that changes only when $\Theta_t^x$ is equal to
$0$ or $\pi$ and $\Theta^x_t$ is always in the interval $[0,\pi]$.
The process $\Theta^x_t$ is constructed in such a way that the
difference $\Theta^x_t - \wt \Theta^x_t$ is constant on every
interval of time on which $\Theta^x_t$ does not hit $0$ or $\pi$.
The analogous reflected process obtained from $\wt \Theta^y_t$
will be denoted $\wh \Theta^y_t$. Let $\tau^\Theta$ be the
smallest $t$ with $\Theta^x_t = \wh \Theta^y_t$. Then we let
$\Theta^y_t = \wh \Theta^y_t$ for $t \leq \tau^\Theta$ and
$\Theta^y_t = \Theta^x_t$ for $t > \tau^\Theta$. We define a
``clock'' by $\sigma(t) = \int_0^t R^{-2}_s ds$. Then $X_t = (R_t,
\Theta^x_{\sigma(t)})$ and $Y_t = (R_t, \Theta^y_{\sigma(t)})$ are
reflected Brownian motions in $D_*$ with normal reflection---one
can prove this using the same ideas as in the discussion of the
skew-product decomposition for 2-dimensional Brownian motion
presented in [IMK]. Moreover, $X$ and $Y$ behave like free
Brownian motions coupled by the mirror coupling as long as they
are both strictly inside $D_*$. The processes will stay together
after the first time they meet. We call $(X,Y)$ a ``mirror
coupling'' of reflected Brownian motions.

The two processes $X$ and $Y$ in the upper half-plane remain at
the same distance from the origin. Suppose now that $D_*$ is an
arbitrary half-plane, and $x$ and $y$ belong to $D_*$. Let $M$ be
the line of symmetry for $x$ and $y$. Then an analogous
construction yields a pair of reflected Brownian motions starting
from $x$ and $y$ such that the distance from $X_t$ to $M \cap \prt
D_*$ is always the same as for $Y_t$. Let $M_t$ be the line of
symmetry for $X_t$ and $Y_t$. Note that $M_t$ may move, but only
in a continuous way, while the point $M_t \cap \prt D_*$ will
never move. We will call $M_t$ the {\it mirror} and the point $H =
M_t \cap \prt D_*$ will be called the {\it hinge}. The absolute
value of the angle between the mirror and the normal vector to
$\prt D_*$ at $H$ can only decrease.

The next level of generality is to consider a mirror coupling of
reflected Brownian motions in a polygonal domain $ D$. For the
first rigorous construction of a mirror coupling in a domain with
piecewise $C^2$-boundary see [AB]. Earlier applications of mirror
couplings in such domains lacked full justification. A technical
problem that prevents us from generalizing the mirror coupling
construction in a half-plane given above to polygons is that it
may occur, with positive probability, that the two processes are
on two different line segments in the boundary of the domain at
the same time (proving this claim does not seem to be trivial; we
omit the proof because it is not needed in this article). Suppose
that $(X,Y)$ is a mirror coupling in a polygonal domain $D$ and
consider an interval $[t_0,t_1]$ such that for every
$t\in[t_0,t_1]$, either $X_t \notin \prt D$ or $Y_t \notin \prt
D$. Let $I$ be the edge of $\prt D$ which is hit first by one of
the particles after time $t_0$. Let $K$ be the straight line
containing $I$. Since the process which hits $I$ does not ``feel''
the shape of $\prt D$ except for the direction of $I$, it follows
that the two processes will remain at the same distance from the
hinge $H_t = M_t \cap K$. The mirror $M_t$ can move but the hinge
$H_t$ will remain constant as long as $I$ remains the side of
$\prt D$ where the reflection takes place. The hinge $H_t$ will
jump when the reflection location moves from $I$ to another edge
of $\prt D$. The hinge $H_t$ may from time to time lie outside
$\prt D$, if $D$ is not convex.

Our arguments will be based in part on the analysis of all
possible movements of the ``mirror'' $M_t$. If $D$ is a polygonal
domain and only one of the processes is on the boundary of $D$ at
time $t_0$, then the possible movements of the mirror on a small
time interval $[t_0, t_0+\Delta t]$ are described in the above
paragraph. We cannot apply the same analysis to the case when both
processes are on the boundary of $D$ at time $t_0$ so we will
provide an alternative approach in Lemma 2.4 below.

With probability one, reflected Brownian motion never visits any
vertices of the union of polygons $\prt D$, so we will assume that
whenever $X_t\in \prt D$ then $X_t$ lies on a single edge of $\prt
D$.

Suppose that $X_{t}\in \prt D$ and let $K_{X,t}$ be the line
containing the edge of $\prt D$ to which $X_{t}$ belongs. We will
be interested only in the case when $M_t$ is not perpendicular to
$K_{X,t}$. Consider any other straight line $I$ intersecting $M_t$
at a single point $x$. If $M_\cdot$ turns around the hinge $H_t =
M_t \cap K_{X,t}$ so that the (smaller) angle between $M_\cdot$
and $K_{X,t}$ increases, i.e., the two lines become ``more
perpendicular,'' the intersection point of $M_\cdot$ and $I$ will
move into one of the half-lines $I\setminus \{x\}$; we will denote
the closure of this half-line $I^{X,t}$. Let $K_{Y,t}$ and
$I^{Y,t}$ be defined in an analogous way relative to $Y$. If for
some $t$ both processes belong to $\prt D$ then the above
definitions can be applied to $I= K_{X,t}$ and $I=K_{Y,t}$, so
$K_{X,t}^{Y,t}$ and $K_{Y,t}^{X,t}$ are well defined.

\bigskip

\noindent{\bf Lemma 2.4}. {\sl Suppose that $(X,Y)$ is a mirror
coupling of reflected Brownian motions in $D$. With probability
one, for every $t\geq 0$ such that $X_t,Y_t \in \prt D$ and $M_t$
is not perpendicular to any of the lines $K_{X,t}$ and $K_{Y,t}$,
there exists $a = a(t)>0$ such that for $s\in[t,t+a]$, we have
$M_s\cap K_{X,s} \in K_{X,t}^{Y,t}$ and $M_s\cap K_{Y,s} \in
K_{Y,t}^{X,t}$.

}

\bigskip
\noindent{\bf Proof}. Suppose that $\dist(X_0,Y_0) = r_0>0$ and
fix an arbitrarily small $r\in(0,r_0)$. Consider a $\delta\in(0,
r/100)$. First we will modify the mirror coupling $(X,Y)$ as
follows.

Let $T_1=\inf\{t\geq 0: X_t,Y_t \in \prt D\}$ and let
$H^X_t=M_t\cap K_{X,t}$. Let $X^{\delta,1}_{T_1}$ be the point in
$A_1\df D\cap \prt B(H^X_{T_1}, \dist(X_{T_1}, H^X_{T_1}))$ whose
distance from $K_{X,t}$ is $\delta \land \diam(A_1)$. Let
$Y^{\delta,1}_{T_1}=Y_{T_1}$ and let $\{(X^{\delta,1}_t,
Y^{\delta,1}_t), t\geq T_1\}$ be a mirror coupling in $D$ starting
from $(X^{\delta,1}_{T_1}, Y^{\delta,1}_{T_1})$ at time $T_1$ but
otherwise independent of $\{(X_t,Y_t), t\in[0, T_1]\}$. Let $T_2 =
\inf\{t\geq T_1: X^{\delta,1}_t, Y^{\delta,1}_t \in \prt D\}$. We
continue the construction by induction. Suppose that
$(X^{\delta,j}_t, Y^{\delta,j}_t)$ and $T_{j+1}$ have been
defined. Then we let $M^j_t$ denote the mirror for
$X^{\delta,j}_t$ and $ Y^{\delta,j}_t$, and $H^{X,j}_t =M^j_t\cap
K_{X^{\delta,j},t}$. We define $X^{\delta,j+1}_{T_{j+1}}$ to be
the point in $A_{j+1} \df D\cap \prt B(H^{X,j}_{T_{j+1}},
\dist(X^{\delta,j}_{T_{j+1}}, H^{X,j}_{T_{j+1}}))$ whose distance
from $K_{X^{\delta,j},T_{j+1}}$ is $\delta\land \diam(A_{j+1})$.
We also let $Y^{\delta,j+1}_{T_{j+1}}=Y^{\delta,j}_{T_{j+1}}$ and
$\{(X^{\delta,j+1}_t, Y^{\delta,j+1}_t), t\geq T_{j+1}\}$ be a
mirror coupling in $D$ starting from $(X^{\delta,j+1}_{T_{j+1}},
Y^{\delta,j+1}_{T_{j+1}})$ at time $T_{j+1}$ but otherwise
independent of $(X^{\delta,k},Y^{\delta,k})$, $k=1,2,\dots, j$. It
is easy to see that $\sup_j T_j \to \infty$ as $\delta\to0$ in
probability, because reflected Brownian motion does not hit
vertices of $D$.

Let $U_r= \inf_j\inf\{t\geq T_j: \dist(X^{\delta,j}_t,
Y^{\delta,j}_t) \leq r\}$. Let $X^\delta_t = X_t$ for $t\in[0,
T_1)$, $X^\delta_t = X^{\delta,j}_t$ for $t\in[T_j\land U_r,
T_{j+1}\land U_r)$, and $X^\delta_t = X^{\delta,k}_t$ for $t\geq
U_r$, where $k$ is such that $U_r \in[T_k, T_{k+1})$. We define
$Y^\delta$ in a similar way. Note that $Y^\delta$ is a reflected
Brownian motion in $D$ so, trivially, $Y^{1/n}$ converge in
distribution to the reflected Brownian motion in $D$ as $n\to
\infty$.

Before time $U_r$, the distance between $X^\delta$ and $Y^\delta$
is bounded below by $r$ so simple geometry shows that the jumps of
$X^\delta$ at times $T_j < U_r$ satisfy the assumptions of Lemma
2.3. That lemma and a localization argument show that $X^{1/n}$
converge in distribution to a reflected Brownian motion. By
passing to a subsequence, if necessary, we see that $(X^{1/n},
Y^{1/n})$ converge in distribution to $(X^*,Y^*)$, where $X^*$ and
$Y^*$ are reflected Brownian motions in $D$. It follows easily
from the definition of $(X^\delta, Y^\delta)$ that $(X^*,Y^*)$ is
a mirror coupling on every interval $[t_0, t_1]$ such that neither
$X^*_t$ nor $Y^*_t$ visit $\prt D$ for $t\in[t_0,t_1]$. By the
uniqueness of the mirror coupling proved in [AB], it follows that
$(X^*,Y^*)$ has the same distribution as $(X,Y)$.

Let $M^\delta_t$, $K_{\delta,X,t}$, etc., be defined relative to
$(X^\delta, Y^\delta)$ in the same way as $M_t$, $K_{X,t}$, etc.,
have been defined for $(X,Y)$. The jumps of $X^\delta$ have been
chosen so that $M^\delta_t\cap K_{\delta,X,t} $ moves in one
direction along $ K_{\delta,X,t} $, and the same holds for
$M^\delta_t\cap K_{\delta,Y,t} $ and $ K_{\delta,Y,t} $, as long
as $X^\delta$ and $Y^\delta$ are reflecting on the same two edges
of $\prt D$. It is not hard to see that this property is preserved
under the passage to the limit in distribution and that it implies
the statement in the lemma.
 \qed

\bigskip

From now on, we will restrict our attention to the domain $D_1 \df
\{(x_1,x_2)\in D: x_1 >0, x_2 > -1\}$. Let $\prt^d D_1 =
\{(x_1,x_2) \in \prt D_1: x_2 =-1\}$, $\prt^n D_1 = \prt D_1
\setminus \prt^d D_1$, $\prt^\ell D_1 = \{(z_1,z_2)\in \prt D_1:
z_1 =0\}$, and $\prt^s D_1 = \prt D_1 \setminus (\prt^d D_1 \cup
\prt^\ell D_1)$.

Consider the restriction of $\vphi$ to $D_1$ normalized so that
$\vphi(x)>0$ in $D_1$. By Lemma 2.2, $\vphi$ is an eigenfunction
for the Laplacian in $D_1$ with the following boundary conditions:

\bigskip

\item{(2.3)} Dirichlet boundary conditions on $\prt^d D_1 $, and
Neumann boundary conditions on $\prt^n D_1 $.

\bigskip

Since $\vphi(x)>0$ for $x\in D_1$, $\vphi$ is the first
eigenfunction in $D_1$ with the boundary conditions (2.3). By
Lemma 2.2, it will suffice to show that $\vphi(0,0)> \sup_{x\in
\prt^s D_1} \vphi(x)$ to prove Theorem 1.1.

For $z\in D_1$, let $\rho(z)$ denote the infimum of lengths of
Jordan arcs contained in $D_1$ and joining $z$ with $\prt^d D_1$.

\bigskip

\noindent{\bf Lemma 2.5}. {\sl Suppose that $x,y\in \ol D_1$,
$x=(x_1,x_2)$, $y=(y_1,y_2)$, and the line of symmetry $M$ for $x$
and $y$ does not intersect $ \prt^\ell D_1$. Assume that one of
the following conditions holds, (i) $\angle(y-x) \in
[-\pi/4,\pi/4]$ and $x_1\leq 2$, or (ii) $\angle(y-x, x-(2,-1))
\in [-3\pi/4, -\pi/4]$ and $x_1 \geq 2$. Then $\vphi(x) \geq
\vphi(y)$.

}

\bigskip
\noindent{\bf Proof}. Let $u(t,z)$ be the heat equation solution
in $D_1$ with the boundary conditions (2.3) and the initial
condition $u(0,z) = 1$ for all $z$. Suppose that $(X,Y)$ is a
mirror coupling of reflected Brownian motions in $ D_1$ with
$(X_0,Y_0) = (x,y)$. The following representation of the heat
equation solution is well known, $u(t,x) = P(T^X_{\prt^d D_1}> t)$
and $u(t,y) = P(T^Y_{\prt^d D_1} > t)$. Suppose that we can show
that $T^Y_{\prt^d D_1} \leq T^X_{\prt^d D_1}$, a.s. Then $u(t,y)
\leq u(t,x)$ for all $t\geq 0$ and the eigenfunction expansion
applied for large $t$ shows that $\vphi(y) \leq \vphi(x)$ (see
Proposition 2.1 and the proof of Theorem 3.3 in [BB1]). Hence, it
will suffice to show that $T^Y_{\prt^d D_1} \leq T^X_{\prt^d
D_1}$, a.s.

Recall the definition of $\rho(z)$ stated before the lemma. It is
enough to show that $\rho(Y_t) \leq \rho(X_t)$ for all $t$. Recall
that $M_t$ denotes the mirror, i.e., the line of symmetry for
$X_t$ and $Y_t$. Suppose that one of the conditions (i) or (ii) in
the statement of the lemma is satisfied by $x$ and $y$. The rules
for the possible movements of $M_t$ described before and in Lemma
2.4 imply that as long as $Y_t \in D_2\df\{(z_1,z_2)\in D_1: z_1
\geq 1/2\}$, the mirror $M_t$ has a tendency to intersect $\prt
D_1$ at angles closer to the right angle than the initial angle,
assuming that the parameter $\eps$ in the definition of $D$ is
small. The proof of the last claim is somewhat tedious but totally
elementary so it is left to the reader. We conclude that $M_t$
cannot turn to the point that $\rho(Y_t)> \rho(X_t)$, as long as
$Y_t \in D_2$.

It remains to analyze possible motions of $M_t$ when $Y_t \notin
D_2$. Then one of the processes may reflect on $\prt^\ell D_1 $
while the other is not too far from $\prt^\ell D_1$. We will show
that $M_t$ never intersects $\prt^\ell D_1$. Note that if $Y_t
\notin D_2$, $\rho(Y_t) \leq \rho(X_t)$ and one of the processes
reflects on $\prt^s D_1$ then the hinge will stay at a fixed point
on $\prt ^s D_1$ and the mirror will move in such a way that its
other intersection point with $\prt D_1$ will not touch $\prt^\ell
D_1$.

Suppose that $Y_t \notin D_2$, $\rho(Y_t) \leq \rho(X_t)$, $M_t$
does not intersect $\prt^\ell D_1$, and one of the processes
(necessarily $X$) reflects on $\prt^\ell D_1$. Then the hinge lies
outside $\ol D_1$. Since both processes $X_t$ and $Y_t$ must be in
$\ol D_1$, the geometry of this domain makes it impossible for
$M_t$ to turn closer to the horizontal direction than $\pi /8$
(actually, this lower bound is closer to $\pi/4$, if $\eps$ is
small). This implies that the relation $\rho(Y_t) \leq \rho(X_t)$
will remain in force if the reflection point belongs to $\prt^\ell
D_1$. Finally, Lemma 2.4 can be used to show that the above
analysis, based on the assumption that only one process at a time
reflects on the boundary, remains valid when we consider the
situation when both processes reflect at the same time.

We conclude that $\rho(Y_t) \leq \rho(X_t)$ for all $t\geq 0$ and
this completes the proof.
 \qed

\bigskip

\noindent{\bf Lemma 2.6}. {\sl Let $a= \sup_{(x_1,x_2)\in D_1,
x_1=1}\vphi(x)$, $\Gamma =\{x\in D_1: \vphi(x) = a\}$, $r_1
=\inf_{(x_1,x_2) \in\Gamma} x_1$, and $r_2 =\sup_{(x_1,x_2)
\in\Gamma} x_1$. Then for small $\eps$ we have $ 1-2\eps\leq r_1
\leq r_2 \leq 1$, and $\inf_{(x_1,x_2)\in \ol D_1, x_1\leq 1/2 }
\vphi(x) \geq \sup_{(x_1,x_2)\in \ol D_1, x_1 \geq r_1} \vphi(x)$.

}

\bigskip
\noindent{\bf Proof}. It follows easily from Lemma 2.5 that
$\angle(\nabla \vphi)\in [3\pi/4, 5\pi/4]$ for $x=(x_1,x_2)\in
D_1$ with $1/4\leq x_1 \leq 3/2$, assuming that $\eps$ is small.
This and simple geometry imply the lemma.
 \qed

\bigskip

Recall the definition of $\rho(x)$ stated before Lemma 2.5.

\bigskip

\noindent{\bf Lemma 2.7}. {\sl Let $\Gamma_1$ and $\Gamma_2$
denote the two connected components of $\prt ^s D_1$. If $x,y \in
\Gamma_1$ and $\rho(x) > \rho(y)$ then $\vphi(x) > \vphi(y)$. A
similar statement holds for $\Gamma_2$.

}

\bigskip
\noindent{\bf Proof}. Suppose that $x,y \in \Gamma_1$ and
$\rho(x)> \rho(y)$. By the proof of Lemma 2.5, $\vphi(x) \geq
\vphi(y)$. In fact, this is all we need to prove Theorem 1.1 but
we will show that the inequality is strict because the proof is
short and easy. Suppose that $\vphi(x) = \vphi(y)$. It is easy to
see that one can find a non-empty open set $A\subset D_1$ such
that for any $z\in A$, the pair $(x,z)$ satisfies the assumptions
of Lemma 2.5, and the same holds for the pair $(z,y)$. By Lemma
2.5, $\vphi(x)\geq \vphi(z) \geq \vphi(y)$. Since we have assumed
that $\vphi(x) = \vphi(y)$, we see that $\vphi(x) = \vphi(z)$ for
all $z\in A$. The remark following Corollary (6.31) in [F] may be
applied to the operator $\Delta+\lambda$ to conclude that  the
eigenfunctions are real analytic and therefore they cannot be
constant on an open set unless they are constant on the whole
domain $D$. This contradiction completes the proof.
 \qed

\bigskip

We will now define a coupling $(X,Y)$ of reflected Brownian
motions in $ D_1$ with $X_0 = (0,0)$ and $Y_0 = (0,\eps)$. The
mechanism of the coupling will change, as time goes on, depending
on the outcome of some events. Let
 $$\eqalign{
 A_1 &= \{(x_1,x_2) \in D_1: 0<x_1 < \eps, x_2 >
 7\eps/10\},\cr
 A_2 &= \{(x_1,x_2) \in \prt A_1: x_1 = \eps, 7\eps/10 \leq x_2\leq
 8\eps/10\},\cr
 A_3 &=  \prt A_1\cap D_1,\cr
 A_4 &= \{(x_1,x_2) \in D_1: 0<x_1 < \eps, -3\eps/10 < x_2 <
 \eps/10\}.
 }$$
Let $\{(X^1_t, Y^1_t), t\geq 0\}$ be a synchronous coupling of
reflected Brownian motions in $D_1$ with $X^1_0 = (0,0)$ and
$Y^1_0 = (0,\eps)$ (see (2.1)-(2.2)). Let
 $$\eqalign{
 S_0&= T^{Y^1}_{ A_3} \land T^{X^1}_{\prt A_4\cap D_1},\cr
 G_0&=\{Y^1_{T^{Y^1}_{ A_3}} \in A_2,
 T^{Y^1}_{ A_3}\leq T^{X^1}_{\prt A_4\cap D_1}\}.
 }$$
We let $(X_t,Y_t) = (X^1_t,Y^1_t)$ for $t\leq S_0$. If $G_0$ does
not occur we let $\{(X_t, Y_t), t\geq S_0\}$ be a mirror coupling
starting from $(X^1_{S_0} , Y^1_{S_0})$, but otherwise independent
from $\{(X_t,Y_t), t\in[0, S_0]\}$. Let $S_{1}^0 =S_2^0=T^{Y^1}_{
A_3}$ and for integer $j\in[0, 2/\eps]$,
 $$\eqalign{
 A_5^j & = \{(x_1,x_2) \in D_1: j \eps < x_1 < (j+2) \eps,
 6\eps/10 < x_2 < 9\eps/10\},\cr
 A_6^j &= \{(x_1,x_2)\in \prt A_5^j: x_1 = (j+2)\eps,
 7\eps/10 < x_2 < 8\eps/10\},\cr
 S_1^{j+1} &= \inf\{t\geq S_2^{j}: Y^1_t\in A_6^j\},\cr
 S_2^{j+1} &= \inf\{t\geq S_2^{j}: Y^1_t\in \prt A_5^j\},\cr
 F_j &= \{S_1^{j+1} \leq S_2^{j+1} \}.
 }$$
Fix some $c_*\in(0,1)$ whose value will be chosen later, let $j_0$
be the integer part of $c_*/\eps$ and $F_* = \bigcap_{0\leq j \leq
j_0} F_j$.

If $G_0$ holds and there exists $j\leq j_0$ such that $F_j$ does
not occur then we let $j_1$ be the smallest $j$ with this
property, $(X_t,Y_t) = (X^1_t,Y^1_t)$ for $t\in[S_0,
S^{j_1+1}_2]$, and $\{(X_t, Y_t), t\geq S^{j_1+1}_2\}$ be a mirror
coupling starting from $(X^1(S^{j_1+1}_2) , Y^1(S^{j_1+1}_2))$,
but otherwise independent of $\{(X_t,Y_t), t\in[0,
S^{j_1+1}_2]\}$. Let
 $$\eqalign{
 A_7 &= \{(x_1,x_2)\in D_1: j_0\eps < x_1 < (j_0+3)\eps,
 x_2 > -5\eps/10\},\cr
 A_{8} &= \{(x_1,x_2)\in\prt A_7 \cap \prt D_1: (j_0+1)\eps \leq
 x_1 \leq (j_0+2)\eps\},\cr
 S_3 &= \inf\{t\geq S^{j_0+1}_2: X^1_t  \in \prt A_7\}, \cr
 G_1 &= \{ X^1_{S_3} \in A_{8}\}.
 }$$
If $G_0$ and $F_*$ hold then we let $(X_t,Y_t) = (X^1_t,Y^1_t)$
for $t\in[S_0, S_3]$. We let $\{(X^2_t, Y^2_t), t\geq S_3\}$ be a
mirror coupling starting from $(X^1_{S_3}, Y^1_{S_3})$, but
otherwise independent of the process $\{(X_t,Y_t), t\in[0,
S_3]\}$. Let
 $$\eqalign{
 A_{9} &= \{(x_1,x_2)\in D_1: (j_0-100)\eps < x_1 <
 (j_0+100)\eps, x_2>0\}, \cr
 A_{10} &= \{(x_1,x_2)\in \prt A_{9}: x_1 = (j_0-100)\eps
 \hbox{  or  } x_1 = (j_0+100)\eps\},\cr
 S_4 &= \inf\{t\geq S_3: X^2_t \in A_{10}\},\cr
 S_5 &= \inf\{t\geq S_3: Y^2_t \in A_{10}\},\cr
 S_6 &= \inf\{t\geq S_3: X^2_t \in \prt A_{9}\},\cr
 S_7 &= \inf\{t\geq S_3: Y^2_t \in \prt A_{9}\},\cr
 S_8 &= S_6\land S_7,\cr
 G_2 &= \{S_4 \leq S_5 \land S_6\} \cup
 \{ S_5 \leq S_4 \land S_7\}.
 }$$
If $G_0$ and $F_*$ hold then we let $(X_t,Y_t) = (X^2_t,Y^2_t)$
for $t\in[S_3,S_8]$. If $G_0$ and $F_*$ hold but $G_2$ does not
occur then we let $(X_t,Y_t) = (X^2_t,Y^2_t)$ for $t\geq S_8$. We
let $\{(X^3_t, Y^3_t), t\geq S_8\}$ be a pair of reflected
Brownian motions in $D_1$ starting from $(X^2_{S_8} , Y^2_{S_8})$,
independent from each other and independent from $\{(X_t,Y_t),
t\in[0, S_8]\}$. Let
 $$\eqalign{
 A_{11} &=\{(x_1,x_2)\in D_1: (j_0-50)\eps < x_1 <
 (j_0+50)\eps \hbox{  or  } x_1<\eps\}, \cr
 A_{12} &= \{(x_1,x_2)\in D_1:  x_1 = 1\},\cr
 A_{13} &= \{(x_1,x_2)\in D_1: x_1 \leq 1/2\}, \cr
 S_9 &= \inf\{t\geq S_8: X^3_t \in  A_{11}\},\cr
 S_{10} &=
 \inf\{t\geq S_8: Y^3_t \in  A_{11}\},\cr
 S_{11} &=
 \inf\{t\geq S_8: Y^3_t \in A_{12}\},\cr
 S_{12} &= S_9\land S_{10},\cr
 G_3 &= \{ S_{11} < S_{12}, X_{S_{11}}\in A_{13} \}.
 }$$
If $G_0\cap F_*\cap G_2$ holds then we let $(X_t,Y_t) =
(X^3_t,Y^3_t)$ for $t\in[ S_8, S_{12}]$ and we let $\{(X_t, Y_t),
t\geq S_{12}\}$ be a mirror coupling starting from $(X^3_{S_{12}},
Y^3_{S_{12}})$ but otherwise independent from $\{(X_t, Y_t),
t\in[0, S_{12}]\}$.

\bigskip

\noindent{\bf Lemma 2.8}. {\sl Let $\Gamma$ be the curve defined
in Lemma 2.6 and let $(X,Y)$ be the coupling of reflected Brownian
motions defined before this lemma. There exist $c_1,\eps_1>0$ such
that for $\eps\in(0,\eps_1)$ we have $P(T^X_\Gamma < T^Y_\Gamma)
\leq e^{-c_1/\eps}$.

}

\bigskip
\noindent{\bf Proof}. The pair $(X,Y)$ is not a mirror coupling
but we can still define the ``mirror'' $M_t$ for $(X_t,Y_t)$ as
the line of symmetry for these processes.

Let $K_t=(K^1_t,K^2_t)$ be that of intersection points of the
mirror $M_t$ with $\prt D_1$ which satisfies $\angle(K_t-X_t,
Y_t-X_t) \geq 0$. First we will show that $K^2_t \geq 0$ for all
$t\leq T^X_\Gamma \land T^Y_\Gamma$, a.s., that is, the ``left''
(looking from $X_t$ towards $Y_t$) intersection point of $M_t$
with $\prt D_1$ cannot cross $\prt^\ell D_1$ below $(0,0)$.

We will start by analyzing possible movements of $(X,Y)$. We will
use the following convention, introduced in Lemma 2.3, to denote
coordinates of processes: $X_t = (\wt X_t, \wh X_t)$, and
similarly for other processes. We will show that $\wt X_t \leq \wt
Y_t$ for $t\in[0, S_0]$. Suppose that there is $t_0\in[0,S_0]$
with $\wt X_{t_0} > \wt Y_{t_0}$ and let $t_1 = \sup\{t < t_0: \wt
X_t\leq \wt Y_t\}$. By continuity of reflected Brownian paths,
$\wt X_{t_1} = \wt Y_{t_1}$. Let $W$ be the Brownian motion
driving $X$ and $Y$, in the sense of (2.1)-(2.2). Since $\wt X_t >
\wt Y_t \geq 0$ for $t\in[t_1,t_0]$, $X$ is not reflecting on this
interval, so $\wt X_{t_0} - \wt X_{t_1} = \wt W_{t_0} - \wt
W_{t_1}$. The horizontal component of the vector of reflection for
$Y$ is non-negative so $\wt Y_{t_0} - \wt Y_{t_1} \geq \wt W_{t_0}
- \wt W_{t_1}$. We see that $\wt Y_{t_0} - \wt Y_{t_1} \geq \wt
X_{t_0} - \wt X_{t_1}$ and this contradicts the facts that $\wt
X_{t_0} > \wt Y_{t_0}$ and $\wt X_{t_1} = \wt Y_{t_1}$. Hence,
$\wt X_t \leq \wt Y_t$ for $t\in[0, S_0]$. This implies that
$K^2_t \geq 0$ for $t\in[0, S_0]$.

Recall the definitions of $j_1$ and $j_0$ from the construction of
$(X,Y)$. On the interval $[S_0, S_2^{j_1+1}]$, processes $X$ and
$Y$ do not hit the boundary of $D_1$, so the mirror is translated
but not rotated and the constraints on the positions of $X$ and
$Y$ are such that it is easy to see that $K^2_t \geq 0$ for
$t\in[S_0, S_2^{j_1+1}]$.

Suppose that $G_0\cap F_*$ holds. Then only $Y$ can be reflecting
on the interval $[S_2^{j_0+1}, S_3]$, so $\wt Y_t - \wt X_t$ is
non-decreasing on $[S_2^{j_0+1}, S_3]$. It follows that $\wt X_t
\leq \wt Y_t$ for $t\in[S_2^{j_0+1}, S_3]$ and $K^2_t \geq 0$ on
this interval.

If $G_0 \cap F_* \cap G_1 \cap G_2$ holds then $\wt X_t \leq \wt
Y_t$ for $t\in[S_3, S_{12}]$, so $K^2_t \geq 0$ on this interval.

Suppose that $K^2_t = 0$ for some $t\leq T^X_\Gamma \land
T^Y_\Gamma$ and let $U=\inf\{t\geq 0: K^2_t=0\}$. The above
analysis covers all cases when $X$ and $Y$ are not mirror-coupled.
In other words, if $U$ exists then $\{(X_t,Y_t), t\geq U\}$ is a
mirror coupling. It is not hard to see that an even stronger
statement holds---for some $U_1 < U$, $\{(X_t,Y_t), t\geq U_1\}$
is a mirror coupling.

Suppose that $K^1_U >2$. This means that $\angle(Y_t-X_t)$ must
have changed its value from $\pi/2$ at time $t=0$, to 0 or $\pi$
at some time $T_0\leq T^X_\Gamma \land T^Y_\Gamma$, and then take
a value less then $-\pi/4$ at time $U$. Such a change of
$\angle(Y_t-X_t)$ between $T_0$ and $U$ is impossible, by the
argument given in the proof of Lemma 2.5. Next suppose that $K_U =
(0,0)$. If $\angle(Y_U-X_U)=\pi/2$ then, by symmetry and
uniqueness of the mirror coupling, $\angle(Y_t-X_t)=\pi/2$ and
$M_t = M_U$ for all $t\in[U, T^X_\Gamma \land T^Y_\Gamma]$. Hence,
in this case, $K^2_t \geq 0$ for all $t\leq T^X_\Gamma \land
T^Y_\Gamma$.

Suppose that $K_U = (0,0)$ and $\angle(Y_U-X_U)<\pi/2$. Note that
at time $U$, at least one of the processes must be on the boundary
of $D_1$ (otherwise the mirror is not moving). In the present
case, geometry shows that $X_U\in \prt D_1$ and $Y_U\notin \prt
D_1$. Hence, for some $U_2 < U$ and all $t\in[U_2,U]$, $Y_t \notin
\prt D_1$. This implies that the only process that can reflect on
$\prt D_1$ on the interval $[U_2,U]$ is $X$. However, such
reflection could only push $K^2_t$ up, so $K^2_t \leq 0$ for
$t\in[U_2,U]$, a contradiction with the definition of $U$.

Now assume that $K_U = (0,0)$ and $\angle(Y_U-X_U)>\pi/2$. The
point $v\df(-1,0)$ lies at the intersection of lines containing
the two line segments comprising $J\df \{(x_1,x_2)\in \prt D_1: 0
< x_1 < 1\}$. Let $Z_t $ be the intersection point of $M_t$ with
$\{(x_1,x_2): x_1=-1\}$. Note that the introductory arguments in
this proof showed not only that $U$ cannot occur when $X$ and $Y$
are not mirror-coupled, but also that $M_t$ passes above $v$
(i.e., $Z_t>0$) for $t\in[0,U_1]$. Since $M_U$ passes below $v$
(i.e., $Z_U<0$), there must be a time $U_3 \in (U_1,U)$ such that
either $M_{U_3}$ is vertical and $\angle(Y_t-X_t) = 0$ or $v\in
M_{U_3}$. In the first case, we have $\angle(Y(T^X_\Gamma \land
T^Y_\Gamma) - X(T^X_\Gamma \land T^Y_\Gamma)) \in [-\pi/4,\pi/4]$,
by the argument given in the proof of Lemma 2.5. In the second
case, let $U_4=\inf\{t\geq U_1: Z_t < 0\}$. By continuity, $v\in
M_{U_4}$. At time $U_4$, at least one of the processes must be on
the boundary. Since $U_4<U$, $X_{U_4} \notin \prt D_1$. On a small
interval $[U_4,U_5]$, only $Y$ can reflect on the boundary of
$\prt D_1$. But this reflection will either leave $Z_t$ unchanged
(if $Y$ is reflecting on $J$) or it will push $Z_t$ up (if $Y$ is
reflecting on $\prt ^\ell D_1$), and this contradicts the
definition of $U_4$. This completes the proof of the claim that
$K^2_t \geq 0$ for all $t\leq T^X_\Gamma \land T^Y_\Gamma$.

Let $C$ be the part of $ \{(x_1,x_2)\in \prt D_1: 0 \leq x_1 \leq
1, x_2>0\}$ that lies to the left of $\Gamma$, and $V=T^X_C$.
Standard estimates show that $P(V\geq T^X_\Gamma)\leq
e^{-c_1/\eps}$ for some $c_1>0$. If $V< T^X_\Gamma$ and $V\geq
T^Y_\Gamma$ then $T^X_\Gamma\geq T^Y_\Gamma$. Suppose that $V<
T^X_\Gamma$ and $V<T^Y_\Gamma$. Then we have three possibilities.
First, $X_V=Y_V$. This implies that $T^X_\Gamma=T^Y_\Gamma$. The
second possibility is that $\angle(Y_V-X_V) >\pi$. This implies
the existence of a time $U_6 < V$ such that $K_{U_6} =(0,0)$,
which is impossible by the first part of the proof. Finally,
suppose that $\angle(Y_V-X_V) <\pi$. Then in fact $\angle(Y_V-X_V)
<\pi/8$, for small $\eps$. This implies that $T^X_\Gamma\geq
T^Y_\Gamma$, by an argument similar to that in the proof of Lemma
2.5, because $\{(X_t,Y_t), t\geq V\}$ is necessarily a mirror
coupling, except for the interval $[S_8, S_{12}]$, where the
processes are independent but well separated. We conclude that
$T^X_\Gamma<T^Y_\Gamma$ only if $\{V\geq T^X_\Gamma\}$ occurs.
Since the probability of this event is bounded by $e^{-c_1/\eps}$,
the lemma follows.
 \qed

\bigskip

\noindent{\bf Lemma 2.9}. {\sl There exist $c_2,\eps_1>0$ such
that for $\eps\in(0,\eps_1)$ and $t\geq 1$ we have $P(t\leq
T^X_\Gamma < T^Y_\Gamma) \leq e^{-c_2 t/\eps^2}$.

}

\bigskip
\noindent{\bf Proof}. Recall the set $C$ from the proof of Lemma
2.8 and let $C_0$ be the part of $D_1$ to the left of $\Gamma$.
Let $Q_j $ be the event that $X$ does not hit $C$ during the time
interval $[j\eps^2,(j+1)\eps^2]$. It is easy to see that
$P(Q_j\mid X_{j\eps^2} \in C_0, {\cal F}_{j\eps^2}) < p_1<1$,
where $p_1$ is independent of $j$. By the Markov property applied
at times $j\eps^2$, $P(\{T^X_\Gamma \geq t\} \cap \bigcap_{j\leq
t/\eps^2} Q_j ) \leq p_1^{t/\eps^2} = e^{-c_2 t/\eps^2}$, for some
$c_2>0$. If one of the events $ Q_j^c$ does happen, an argument
similar to that in the proof of Lemma 2.8 shows that
$T^X_\Gamma\geq T^Y_\Gamma$.
 \qed

\bigskip

\noindent{\bf Lemma 2.10}. {\sl Let $C_1=\{(x_1,x_2)\in D_1: x_1
\leq 1/2\}$. For any $c_3>0$ there exist $c_*\in (0,1)$ (used in
the construction of the coupling $(X,Y)$) and $\eps_1>0$ such that
for $\eps\in(0,\eps_1)$ we have $P(T^Y_\Gamma < T^X_\Gamma,
X_{T^Y_\Gamma}\in C_1) \geq e^{-c_3/\eps}$.

}

\bigskip
\noindent{\bf Proof}. We will argue that for some constants $c_j$
independent of $\eps$ we have the following bounds for the
probabilities of events defined in the construction of the
coupling $(X,Y)$,
 $$\eqalignno{
 P(G_0) &\geq c_4,&(2.4)\cr
 P\left(F_j\mid G_0 \cap \bigcap_{k<j} F_k\right) &\geq c_5,
 \qquad 0\leq j \leq j_0,&(2.5)\cr
 P\left(G_1\mid G_0 \cap \bigcap_{0\leq j \leq j_0} F_k\right)
 &\geq c_6, &(2.6)\cr
 P\left(G_2\mid G_0 \cap G_1\cap \bigcap_{0\leq j \leq j_0} F_k\right)
 &\geq c_7\eps,&(2.7)\cr
 P\left(G_3\mid G_0 \cap G_1\cap G_2\cap
 \bigcap_{0\leq j \leq j_0} F_k\right) &\geq c_8 \eps^2.&(2.8)
 }$$

Let $W$ be the Brownian motion driving $(X,Y)$, in the sense of
(2.1)-(2.2), on the interval $[0, S_0]$. Recall the notation
$W=(\wt W, \wh W)$. By the support theorem (see Theorem I (6.6) in
[B]) for the planar Brownian motion and scaling, the following
event $G_0^*$ has probability greater than $c_4>0$, independent of
$\eps$.

\item{($G_0^*$)} The Brownian motion $W$ goes from $(0,0)$ to
$B((0, -0.25\eps), 0.01\eps)$ before touching the boundary of
$\{(x_1,x_2)\in \R^2: |x_1| < 0.02\eps, -0.26\eps < x_2 < 0.01
\eps\}$ in less than $\eps^2$ units of time, and then goes to
$B((2\eps, -0.25\eps), 0.01\eps)$ without hitting the boundary of
$\{(x_1,x_2)\in\R^2: -0.01\eps < x_1 < 3\eps, -0.26 \eps < x_2 <
-0.24 \eps\}$, in another time interval of $\eps^2$ units or less.

\medskip Let $T_*$ be the time needed to complete the movements
described in $G^*_0$. We will argue that if $G^*_0$ holds then so
does $G_0$. Note that $\wh X = \wh W$ on $[0, T_*]$. We have
already shown that $\wt X_t \leq \wt Y_t$ for $t\in [0, S_0]$ in
the proof of Lemma 2.8, so it remains to show that $S_0\leq T_*$
and $Y_{T^{Y}_{ A_3}} \in A_2$. Since $\wh W_t \leq  0.01 \eps$
for $t\in [0, T_*]$, we have $\wh Y_t \in [\wh W_t - 0.01 \eps,
\wh W_t]$ for $t\in [0, T_*]$. The reflection vector for $Y$ is
either horizontal or pointing down, at an angle not greater than $
\eps$ with the vertical. Since $\wt W_t \geq -0.02 \eps$ for $t\in
[0, T_*]$, this implies that $\wt Y_t \in [\wt W_t, \wt W_t+ 0.01
\eps^2 + 0.02\eps]$ for $t\in [0, T_*]$. Now easy geometry shows
that $X$ and $Y$ are transforms of $W$ that satisfy $G_0$. This
proves (2.4).

The support theorem for the planar Brownian motion (i.e., without
reflection) easily yields (2.5) and (2.6).

If $ G_0 \cap G_1\cap \bigcap_{0\leq j \leq j_0} F_k$ holds then
$\dist(X_{S_3}, Y_{S_3}) \geq c_9 \eps^2$ because $X$ is located
to the left of $Y$ (in the sense of the first coordinate) at time
$S_2^{j_0+1}$. The event $G_2$ will occur if $Y$ moves above the
horizontal axis, about $100 \eps$ units to the right before moving
$c_9 \eps^2$ units to the left. By the ``gambler's ruin''
estimate, we obtain (2.7).

Recall the point $v\df(-1,0)$ that lies at the intersection of
lines containing the two line segments comprising $ \{(x_1,x_2)\in
\prt D_1: 0 < x_1 < 1\}$. Let $R^X_t = \dist (X_t,v)$ and $R^Y_t =
\dist (Y_t,v)$. As long as $X$ and $Y$ stay inside $\{(x_1,x_2)\in
D_1: \eps \leq x_1 \leq 1-\eps\}$, $R^X_t$ and $R^Y_t$ are
2-dimensional Bessel processes because the reflection has no
effect on the distance of $X $ or $Y$ from $v$. It is standard to
show the the 2-dimensional Bessel process $R^Y_t$ starting about
$c_{10} + 100\eps$ units from 0 will reach the value $1+\eps$ at a
time $T_0\in(1/2,1)$, before hitting the level $c_{10} + 50\eps$,
with probability exceeding $c_{11} \eps$. The other 2-dimensional
Bessel process, $R^X$, starting about $c_{10} - 100\eps$ units
from 0, will stay in the interval $(1/4,1/2)$ during the time
interval $(1/2,1)$, before hitting levels $c_{10} -50\eps$ or
$\eps$, with probability exceeding $c_{12} \eps$. By independence
of $R^X$ and $R^Y$ on $[S_8, S_{12}]$, we obtain (2.8).

Recall that $j_0$ is the integer part of $c_*/\eps$ and $F_* =
\bigcap_{0\leq j \leq j_0} F_j$. It follows from (2.4)-(2.8) and
the repeated application of the strong Markov property that
 $$P(G_0 \cap F_* \cap G_1 \cap G_2 \cap G_3) \geq c_{13} \eps^3
 c_5^{c_*/\eps}.\eqno(2.9)$$
For any fixed $c_3>0$, the right hand side of (2.9) is greater
than $e^{-c_3/\eps}$ if $c_*$ and $\eps$ are sufficiently small.
If the event in (2.9) occurs then $T^Y_\Gamma < T^X_\Gamma$ and $
X_{T^Y_\Gamma}\in C_1$, so the lemma follows.
 \qed

\bigskip

\noindent{\bf Proof of Theorem 1.1}. Lemmas 2.2 and 2.7 show that
it will suffice to prove that $\vphi(0,0) > \vphi(0,-\eps)\lor
\vphi(0,\eps)$. We will only show that $\vphi(0,0) >
\vphi(0,\eps)$ because the claim that $\vphi(0,0) >
\vphi(0,-\eps)$ can be proved in a completely analogous way, by
symmetry.

Let $(X,Y)$ be the coupling constructed before Lemma 2.6, $S=
T^X_\Gamma \land T^Y_\Gamma$ and let $U$ be the coupling time,
i.e., $U=\inf\{t\geq 0: X_t = Y_t\}$. Let $u(t,x) = \vphi(x)
e^{\lambda t}$ and note that $u$ is a solution to the heat
equation with the Neumann boundary conditions on $\prt^n D_1$ and
Dirichlet boundary conditions on $\prt^d D_1$. Since $X$ and $Y$
are reflected Brownian motions in $D_1$, we have the following
probabilistic representation of $u$, for bounded stopping times
$T\leq T^X_{\prt^d D_1}$,
 $$u(0,(0,0)) = E u(X_T,T) = E\vphi(X_T) e^{\lambda T}.$$
Let $D_2 =\{(x_1,x_2)\in D_1: x_1 \leq 1\}$ and let $\mu$ be the
first eigenvalue for the Laplacian in $D_2$ with Neumann boundary
conditions on $\prt^n D_1 \cap \prt D_2$ and Dirichlet boundary
conditions elsewhere. It is easy to see that $\mu> \lambda$. We
have
 $$P(S>t) \leq P(T^X_{\prt D_2\cap D_1}>t)\leq c_1 e^{-\mu t}.$$
This and the fact that the eigenfunction $\vphi$ is bounded
(because $D$ is Lipschitz) imply that random variables
$\vphi(X_{S\land n}) e^{\lambda (S\land n)}$ are dominated by a
random variable with an exponential tail. Hence, we can use the
fact that $u(0,(0,0)) = E\vphi(X_{S\land n}) e^{\lambda (S\land
n)}$ and the dominated convergence theorem to prove that
$u(0,(0,0)) = E\vphi(X_S) e^{\lambda S}$. Similarly,
$u(0,(0,\eps)) = E\vphi(Y_S) e^{\lambda S}$. We have
 $$\eqalign{
 u(0,(0,0)&) - u(0,(0,\eps))
 = E\vphi(X_S) e^{\lambda S} - E\vphi(Y_S) e^{\lambda S}\cr
 &= E (\vphi(X_S) e^{\lambda S} - \vphi(Y_S) e^{\lambda S})
 \bone_{\{S>U\}}\cr
 & = E (\vphi(X_S) - \vphi(Y_S)) e^{\lambda S}
 \bone_{\{T^X_\Gamma < T^Y_\Gamma\}}
 +  E (\vphi(X_S) - \vphi(Y_S)) e^{\lambda S}
 \bone_{\{T^X_\Gamma > T^Y_\Gamma\}}.
 }$$

Let $a_1 = \sup_{ x=(x_1,x_2)\in D_1, x_1\geq
1-2\eps,y\in\Gamma}\vphi(x) -\vphi(y)$ and recall from Lemma 2.6
that $a_1 \leq a_2 \df \inf_{ x=(x_1,x_2)\in D_1, x_1\leq 1/2,
y\in\Gamma}\vphi(x) -\vphi(y)$. By the proof of Lemma 2.8,
$\vphi(X_S) - \vphi(Y_S) \geq -a_1$ if $T^X_\Gamma < T^Y_\Gamma$.
By Lemmas 2.1, 2.8 and 2.9,
 $$\eqalign{
 E &(\vphi(X_S) - \vphi(Y_S)) e^{\lambda S}
 \bone_{\{T^X_\Gamma < T^Y_\Gamma\}}\cr
 &=E (\vphi(X_S) - \vphi(Y_S)) e^{\lambda S}
 \bone_{\{T^X_\Gamma < T^Y_\Gamma, S\leq 1\}}\cr
 &\qquad + \sum_{k\geq 1} E (\vphi(X_S) - \vphi(Y_S)) e^{\lambda S}
 \bone_{\{T^X_\Gamma < T^Y_\Gamma, S\in[k,k+1]\}}\cr
 &\geq -a_1 e^\lambda e^{-c_1/\eps} -
 \sum_{k\geq 1} a_1 e^{\lambda (k+1)} e^{-c_2 k /\eps^2}.
 }$$
For some $c_4$ and small $\eps$, this is greater than $- a_1 c_4
e^{-c_1/\eps}$. According to Lemma 2.10, one can choose $c_*$ and
$ c_3$ such that for small $\eps$,
 $$
 E (\vphi(X_S) - \vphi(Y_S)) e^{\lambda S}
 \bone_{\{T^X_\Gamma > T^Y_\Gamma\}}\geq a_2 e^{-c_3/\eps}
 > a_1 c_4 e^{-c_1/\eps}.
 $$
Hence $u(0,(0,0)) - u(0,(0,\eps)) >0$ and $\vphi(0,0) >
\vphi(0,\eps)$.
 \qed

\vfill\eject

\vskip1truein \centerline{REFERENCES}
\bigskip

\item{[A]} R. Atar (2001). Invariant wedges for a two-point
reflecting Brownian motion and the ``hot spots'' problem. {\it
Elect. J. of Probab.} 6, paper 18, 1--19.

\item{[AB]} R. Atar and K. Burdzy (2004) On Neumann eigenfunctions
in lip domains. {\it  J. Amer. Math. Soc. \bf 17}, 243--265.

\item{[BB1]} R. Ba\~nuelos and K. Burdzy (1999). On the ``hot
spots'' conjecture of J. Rauch. {\it J. Funct. Anal.} 164, 1--33.

\item{[B]} R.F. Bass, {\it Probabilistic Techniques in Analysis},
Springer, New York, 1995.

\item{[BB2]} R. Bass and K. Burdzy (2000). Fiber Brownian motion
and the `hot spots' problem {\it Duke Math. J.} 105, 25--58.

\item{[BK]} K. Burdzy and W. Kendall (2000) Efficient Markovian
couplings: examples and counterexamples {\it Ann. Appl. Probab.
\bf 10} (2000) 362--409.

\item{[BW]} K. Burdzy and W. Werner (1999). A counterexample to
the "hot spots" conjecture {\it Ann. Math.} 149, 309--317.

\item{[CH]} R.~Courant and D.~Hilbert, {\it Methods of
Mathematical Physics}, Interscience Publishers, New York, 1953.

\item{[F]} G.B.~Folland, {\it Introduction to Partial Differential
Equations}, Princeton Univ. Press, Princeton, 1976.

\item{[IMK]} K.~It\^o and H.P.~McKean (1974) {\it Diffusion
Processes and Their Sample Paths}, Springer, Berlin.

\item{[JN]} D. Jerison and N. Nadirashvili (2000) The ``hot
spots'' conjecture for domains with two axes of symmetry. {\it  J.
Amer. Math. Soc. \bf 13}, 741--772.

\item{[K]} B.~Kawohl, {\it Rearrangements and Convexity of Level
Sets in PDE}, Lecture Notes in Mathematics 1150, Springer, Berlin,
1985.

\item{[LS]} P.~L. Lions and A.~S. Sznitman (1984) Stochastic
differential equations with reflecting boundary conditions. {\it
Comm. Pure Appl. Math.} {\bf 37}, 511-537.

\item{[P]}  M. Pascu (2002) Scaling coupling of reflecting
Brownian motions and the hot spots problem. {\it Trans. Amer.
Math. Soc. \bf 354}, 4681--4702.

\bigskip

Department of Mathematics, Box 354350, University of Washington,
Seattle, WA 98115-4350, burdzy@math.washington.edu

\bye